\RequirePackage{fix-cm}
\documentclass[smallextended]{svjour3}       
\smartqed  
\usepackage{graphicx}
\usepackage[utf8]{inputenc}
\usepackage{acronym}
\usepackage{array}
\usepackage{mdwmath}
\usepackage{mdwtab}
\usepackage{fixltx2e}
\usepackage{amsmath}
\usepackage{amssymb}
\usepackage{pgfplots}
\usepackage{tikz}
\usepackage{ifthen}
\usepackage{pdfcomment}
\usepackage[Gray,squaren,thinqspace,thinspace]{SIunits}
\usepackage{lipsum}
\usepackage{caption}
\usepackage{subcaption}
\usepackage{todonotes}
\usepackage{enumitem}
\usetikzlibrary{calc}
\newcommand {\XYcoord}[2]{ 
	\draw[thick, ->] (0,0) -- +(#1,0) node[anchor=north] {$x$};
	\draw[thick, ->] (0,0) -- +(0,#2) node[anchor=east]{$y$};	
	
}
\newcommand{\dd}[2]{\frac{\mbox{d}{#1}}{\mbox{d}{#2}}}                     
\newcommand{\pd}[2]{\frac{\partial{#1}}{\partial{#2}}}                     
\newcommand{\ud}[1]{\,\mbox{d}{#1}}                                        
\newcommand{\udV}{\ud{{V}}}                                               

\newcommand{\pdt}[1]{\pd{#1}{{t}}}
\newcommand{\fdt}[1]{\frac{\mathrm{d}{#1}}{\mathrm{d}t}}                   
\newcommand{\mvc}[1]{\mathbf{#1}}
\newcommand{\MVP}{\mvc{A}}
\renewcommand{\dd}{\boldsymbol{\delta}}
\newcommand{\pp}{{\bf p}}
\newcommand{\ppexpv} {\pp}

\newcommand{\ppup}{p^{\mathrm{u}}}
\newcommand{\pplow}{p^{\mathrm{l}}}

\newcommand{\dlow}{\delta^{\mathrm{l}}}
\newcommand{\dup}{\delta^{\mathrm{u}}}
\newcommand{\slack}{\mathbf{\xi}}

\newcommand{\upk}{^{(k)}}
\newcommand{\qpt}{\text{ .}}             
\newcommand{\qcm}{\text{ ,}}             
\newcommand{\qsc}{\text{ ;}}             

\begin{document}
	
\title{Robust Shape Optimization of Electric Devices Based on Deterministic Optimization Methods and Finite Element Analysis With Affine Decomposition and Design Elements
}
\titlerunning{Robust shape optimization with affine decomposition and design elements}        
	
\author{Ion Gabriel Ion \and Zeger Bontinck \and Dimitrios Loukrezis \and Ulrich R\"omer \and Oliver Lass \and Stefan Ulbrich \and Sebastian Sch\"ops \and Herbert De~Gersem}
\authorrunning{I.G. Ion, Z. Bontinck, et. al.} 
	
\institute{
	I.G. Ion, Z. Bontinck, D. Loukrezis, U. R\"omer, O. Lass, S. Ulbrich, S. Sch\"ops, H. De~Gersem \at
	Graduate School Computational Engineering, Technische Universit\"at Darmstadt, Dolivostra{\ss}e 15, 64293 Darmstadt, Germany;
	\email{degersem@ce.tu-darmstadt.de}
\and
	I.G. Ion, Z. Bontinck, D. Loukrezis, S. Sch\"ops and H. De~Gersem \at
	Institut f\"ur Theorie Elektromagnetischer Felder, Technische Universit\"at Darmstadt, Schlossgartenstra{\ss}e 8, 64289 Darmstadt, Germany;
	\email{degersem@temf.tu-darmstadt.de}
\and
	U. R\"omer \at
	Institut f\"ur Dynamik und Schwingungen, Technische Universit\"at Braunschweig, Schleinitzstra{\ss}e 20, 38106 Braunschweig, Germany;
	\email{u.roemer@tu-braunschweig.de}
\and
	O. Lass, S. Ulbrich \at
	Chair of Nonlinear Optimization, Department of Mathematics, Technische Universit\"at Darmstadt, Dolivostra{\ss}e 15, 64293 Darmstadt, Germany;
	\email{ulbrich@mathematik.tu-darmstadt.de}
\and
	H. De~Gersem \at
	Wave Propagation and Signal Processing Research Group, KU Leuven -- Kulak, Etienne Sabbelaan 53, 8500 Kortrijk, Belgium;
	\email{herbert.degersem@kuleuven-kulak.be}
}
	
\date{Received: date / Accepted: date}
	
\maketitle
	
\begin{abstract}
In this paper, gradient-based optimization methods are combined with finite-element modeling for improving electric devices. Geometric design parameters are considered by affine decomposition of the geometry or by the design element approach, both of which avoid remeshing. Furthermore, it is shown how to robustify the optimization procedure, i.e., how to deal with uncertainties on the design parameters. The overall procedure is illustrated by an academic example and by the example of a permanent-magnet synchronous machine. The examples show the advantages of deterministic optimization compared to standard and popular stochastic optimization procedures such as, e.g., particle swarm optimization.
\end{abstract}

\section{Introduction}

In almost all electric design procedures, numerical optimization is employed as one of the last design steps in order to optimize the device's performance and efficiency, to minimize its weight and size and to save on material and manufacturing costs. Often, the quality of this optimization step indirectly determines the success of the product and, hence, the market position of the company. The reliability, accuracy and computational cost of the numerical optimization procedure becomes in itself a subject of competition. This paper illustrates that shape optimization can be improved substantially when finite element (FE) analysis procedures are equipped with affine decomposition and design elements, such that well-performing deterministic optimization methods become applicable.

Impressive technical improvements have been achieved by numerical optimization on the basis of magnetic equivalent circuits or 2D and 3D FE models. All have led to highly optimized designs, e.g., for permanent-magnet synchronous machines (PMSMs) in automotive applications.  Since three decades, FE-based optimization has been addressed in several text books (e.g. \cite{Di-Barba_2010aa}) and hundreds of journal articles (see e.g. \cite{Duan_2013aa} and the references therein). Although originally gradient-based methods were preferred (see, e.g., \cite{Russenschuck_1990aa,Weeber_1992aa,Takorabet_1997aa}), already for more than two decades, stochastic algorithms are more popular, see, e.g., \cite{Hameyer_1993aa} and \cite{Lok_2017aa}. The majority of proposed procedures opt for \emph{stochastic} or \emph{population-based} optimization methods, such as genetic algorithms and particle swarm optimization (e.g. \cite{Ma_2015aa}), because they allow to use FE solvers as a black-box, they can easily consider geometric parameters, their parallelization is straightforward and they are more likely to find the global optimum. Stochastic algorithms have been used for robust optimization, have been applied together with surrogate modeling and have been extended to multi-objective optimization problems \cite{Baumgartner_2004aa,Ho_2005aa}. In particular for PMSMs, optimization with stochastic methods became the method of choice \cite{Cassimere_2009aa,Bash_2011aa,Sizov_2013aa}.

The trend toward stochastic optimization combined with FE analysis continues without restraint, as is illustrated by the number of according contributions at recent conferences. This paper partially counteracts this tendency by turning back to deterministic optimization algorithms. These are known to converge faster than stochastic optimization methods, albeit possibly to a local optimum. Moreover, the analysis of gradient based methods is more mature, allowing for a rigorous control of mesh discretization errors, for instance. The main drawback of many deterministic methods is, however, the necessity to provide derivatives, which is particularly cumbersome when optimization according to geometric parameters is pursued. This drawback is here addressed explicitly and is alleviated by affine decomposition of the geometry or by the design element approach. The overall deterministic optimization routine is shown to outperform the most popular stochastic algorithms by factors. Moreover, the optimization method will be robustified to include uncertainties on the design parameters.

The paper is structured as follows: Section~\ref{sec:optimization} recalls the basics of mathematical optimization. It clearly distinguishes between deterministic methods (subsection~\ref{subsec:deterministic}) and particle swarm optimization as a relevant representative of stochastic methods (subsection~\ref{subsec:particle_swarm_optimization}). Furthermore, an extension to robust optimization is discussed in subsection~\ref{subsec:robust_optimization}. Section~\ref{sec:model} deals with FE analysis of magnetodynamic fields. The core parts of the paper are subsection~\ref{subsec:affine_decomposition} about affine decomposition and subsection~\ref{subsec:design_element} about design elements, both facilitating and improving the calculation of derivatives with respect to geometric parameters. The superior performance of gradient-type deterministic optimization is illustrated for a benchmark example in section~\ref{sec:example1} and for a PMSM in section~\ref{sec:example2}. Conclusions are formulated in section~\ref{sec:conclusions}.

\section{Constrained Optimization}
\label{sec:optimization}

\subsection{Constrained optimization problem}
\label{subsec:optimization_problem}

The optimization is carried out with respect to $I$ \emph{design parameters} $\pp=(p_1,p_2,\ldots,p_I)$ belonging to the \emph{admissible set} $\mathcal{P}_{\mathrm{ad}} =\{\pp\in\mathbb{R}^I|G_m(\pp)\leq0,m=1,\ldots,M\}$, where $G_m(\pp)$ denote the constraints. 
The design parameters can be any continuous variables, e.g., material constants, excitation parameters and geometric sizes or positions. The constraints limit the admissible range of these parameters, e.g., to preserve the topology of the geometry or to set physical and operational constraints. Discrete design parameters are not considered in this work, although many methods apply, e.g., as part of a branch-and-bound technique, to mixed-integer optimization problems as well \cite{Hemker_2008aa}.

The optimization goal is represented by the \emph{objective function} $J(\pp)$ returning a scalar value for every set of design parameters. Relevant quantities are, e.g., force, torque, current, efficiency, weight, temperature or a combination thereof. When $Q$ objective functions $J_n(\pp)$, $n=1,\ldots,N$ are relevant, a possible approach is to combine them with user-defined weight factors $\alpha_n$ into a single cost function $J(\pp)=\sum_{n=1}^N \alpha_n J_n(\pp)$. The optimization problem then reads
\begin{subequations}\label{eq:optprb}
	\begin{alignat}{2}
	& \underset{\pp\in \mathbb R^I}{\text{minimize}}	& \quad & J(\pp) \qcm \\
	& \text{subject to}	& \quad & G_m(\pp)       \leq 0 \text{ , } m=1,\ldots,M \qpt
	\end{alignat}
\end{subequations}
In this work, the evaluation of $G_m(\pp)$ and/or $J(\pp)$ involves a FE analysis of the device. Hence, the computational performance of the overall approach is heavily determined by the number of FE-solver calls.

\subsection{Optimization methods}
\label{subsec:optimization_methods}

The selection of a particular optimization method consists of four, essentially independent choices (see also Table~3 in \cite{Graeb_2007aa}).
\begin{itemize}
\item Problem (\ref{eq:optprb}) considers a single optimization goal. For a \emph{multi-objective} optimization problem, a Pareto front is calculated such that the relative importance of the optimization goals can be fixed in a later design stage \cite{Di-Barba_2010aa,Brisset_2008aa}. This paper does not further consider multi-objective optimization. Nonetheless, the developed techniques are applicable to multi-objective optimization as well.
\item Especially when the evaluation of the objective function is computationally expensive, it is recommended to carry out the optimization method on the basis of a \emph{surrogate} model (\emph{indirect} optimization methods). Such a simplified model can be obtained by expert knowledge on the application \cite{Wrobel_2003aa}, by design space reduction~\cite{Gillon_2000aa}, by a response surface methodology~\cite{Gillon_2000aa} or by space mapping~\cite{Koziel_2013aa} or manifold mapping~\cite{Echeverria_2006aa}. Here, a \emph{direct} optimization procedure is used. All ideas presented here can, however, be used in combination with indirect optimization approaches as well.
\item The result from a \emph{nominal} optimization is a set of optimized design parameters leading to an optimum of the objective function. The optimum may, however, become irrelevant when it is highly sensitive to uncertainties in the design parameters. One speaks about \emph{robust} optimization when the optimization is carried out taking such uncertainties into account. In this paper, both nominal and robust optimization methods are considered. An approach for \emph{robustification} is discussed in subsection~\ref{subsec:robust_optimization}.
\item Two families of basic optimization methods exist: \emph{deterministic} and \emph{stochastic} methods. Among the stochastic methods, \emph{genetic algorithms} \cite{Man_2012aa}, \emph{differential evolution} \cite{Neri_2010aa} and \emph{particle swarm optimization} (PSO) \cite{Kennedy_1995aa} are well known.
\end{itemize}
This paper motivates the use of a gradient-based deterministic method for nominal and robust optimization and compares it with a standard particle swarm technique.

\subsection{Gradient-based deterministic method}
\label{subsec:deterministic}

This work proposes to solve (\ref{eq:optprb}) by standard Sequential Quadratric Programming (SQP) with damped Broyden-Fletcher-Goldfarb-Shanno (BFGS) updates for the Hessian approximation \cite{Nocedal_2006aa,Hinze_2009aa}. This method establishes locally a second order convergence, which means that
\begin{equation}
	|J(\pp_{k+1}) - J(\pp_\text{opt})| \leq C|J(\pp_k) - J(\pp_\text{opt})|^2
\end{equation}
for $C>0$ and $k$ the iteration step, which should be sufficiently large. The method, however, requires knowledge about the sensitivities of the objective function with respect to the design parameters, i.e., $\nabla_\pp J(\pp)$ or, alternatively, a locally quadratic approximation of the objective function \cite{Powell_2002aa}. Many FE solution and post-processing routines do not provide this information, especially when geometric design parameters are involved. Therefore, one is tempted to approximate the sensitivities by finite differences as in, e.g., \cite{Russenschuck_1990aa}.
This is, however, known to be particularly cumbersome because of the limited accuracy of the finite differences \cite{Weeber_1992aa}. Even when relying on gradient-free deterministic methods (e.g., \cite{Powell_2002aa,Powell_2009aa}), artifacts caused by FE analysis may hamper the convergence of the optimization routines. Eventually, as apparently the only option, deterministic optimization algorithms are abandoned in favor of stochastic approaches. This paper, however, sticks to gradient-based deterministic methods by complementing the FE simulation procedure with sensitivity information. The problems caused  by the presence of geometric parameters are alleviated by introducing affine decomposition (see subsection~\ref{subsec:affine_decomposition}) or, alternatively, design elements (see subsection~\ref{subsec:design_element}) to the FE procedure.

\subsection{Particle swarm optimization}
\label{subsec:particle_swarm_optimization}

Particle swarm optimization (PSO) \cite{Kennedy_1995aa} belongs to the broad class of stochastic algorithms and is particularly popular for electric machines, see e.g. \cite{Ma_2015aa,Baumgartner_2004aa,Ho_2005aa,Cassimere_2009aa,Bash_2011aa}. In PSO, a set of $Q$ particles indicated by $q=1,\ldots,Q$, moves through the admissible set in the design space in search of an optimum. At each iteration step $k$, the algorithm evaluates the objective function $J(\pp)$ in every particle position $\pp_{k,q}$. The newly obtained values are compared to the previous best values in the individual particle histories and the best value of the entire swarm. The corresponding best sets are denoted by $\hat{\pp}_q$ and $\hat{\pp}_\text{swarm}$ respectively. The velocities of the particles are updated according to
\begin{equation}
	 \mvc{v}_q \leftarrow \underbrace{\omega_0\mvc{v}_q}_\text{1)} + \underbrace{\omega_1 \mathbf{N}_1 (\hat{\pp}_q-\pp_{k,q})}_\text{2)} + \underbrace{\omega_2 \mathbf{N}_2 (\hat{\pp}_\text{swarm}-\pp_{k,q})}_\text{3)} \text{ , }
\end{equation}
where $\omega_0$, $\omega_1$ and $\omega_2$ are swarm characteristic constants and $\mathbf{N}_1$ and $\mathbf{N}_2$ are two random diagonal matrices with elements in $[0,1]$ generated independently and uniformly for each particle at every step, representing the \emph{free will} of the swarm. The components of the velocity update are:	
\begin{enumerate}
	 	\item maintain a part of the current velocity;
	 	\item head towards the particle's best found point ($\hat{\pp}_q$);
	 	\item head towards the swarm's best found point ($\hat{\pp}_\text{swarm}$).
\end{enumerate}
If at some iteration there is a particle that leaves the admissible set, its position is projected on the boundary of the admissible set. Initially, all particles are randomly and uniformly distributed in the admissible set and the initial velocities are set to $0$. The particle swarm is a gradient-free method and works for non-smooth functions as well. The iteration ends when a maximum number of iterations is reached, or when the majority of the particles are close enough to the best point $\hat{\pp}_\text{swarm}$, i.e.,
\begin{equation}
	 \frac{1}{Q}\sum\limits_{q=1}^{Q}\|\hat{\pp}_\text{swarm}-\pp_{k,p}\|_2 < \epsilon \qcm
\end{equation}
with a user-defined tolerance $\epsilon$, or if there is no further change in the global best point $\hat{\pp}_\text{swarm}$ over $N_\text{stall}$ consecutive iterations.

\subsection{Robust optimization}
\label{subsec:robust_optimization}

In a \emph{nominal} optimization, one is looking for the minimum value of an objective function. However, during manufacturing small deviations can occur on the parameters. As a consequence, the optimal solution may become suboptimal in reality. \emph{Robust} optimization searches for an optimum that is not too much affected by the expected parameter deviations \cite{Yoon_1999aa,Omekanda_2006aa}.

One possibility is to optimize such that the worst-case scenario within a stochastic set of possibilities around the optimal design parameters is the best possible.
The robust counterpart of (\ref{eq:optprb}) adopting a worst-case scenario is
\begin{subequations}\label{eq:opt_robust}
	\begin{alignat}{2}
		& \underset{\pp\in \mathbb R^I}{\text{minimize}}	& \quad & \max_{\dd\in U} J(\ppexpv+\dd) \qcm \\
		& \text{subject to}	& \quad & \max_{\dd\in U}G_m(\ppexpv+\dd) \leq 0 \text{ , } m=1,\ldots,M \qpt
	\end{alignat}
\end{subequations}
Here, the uncertainty set for the deviations $\dd$ is defined by
\begin{eqnarray}
U &:=& \{\dd\in \mathbb{R}^n\,|\,\dlow_i\le \dd_i \le \dup_i,\,i=1,\ldots,n\} \nonumber\\ 
&=&\{\dd\in \mathbb{R}^n\,|\,\|\mathbf{D}^{-1}\dd\|_\infty \le 1\} \qcm
\label{eq:unc_set}
\end{eqnarray}
where $\mathbf{D}$ is a scaling matrix and where $\delta_i^l=-\delta_i^u$.

The nested optimization problem formulated by (\ref{eq:opt_robust}) is hard to solve. A numerically feasible optimization problem is obtained by approximating the $\max$ problem, i.e., by applying a first order Taylor approximation of the objective function and the constraints with respect to $\pp$ \cite{Diehl_2006aa}, i.e.,
\begin{align}
\label{eq:RO_lin}
J(\ppexpv+\dd) &\approx  J(\ppexpv) +\nabla_{\ppexpv}J(\ppexpv)\cdot\dd \qsc\\
G_m(\ppexpv+\dd) &\approx  G_m(\ppexpv) +\nabla_{\ppexpv} G_m(\ppexpv)\cdot\dd \qcm
\end{align}
for $m = 1,\ldots,M$. Inserting this approximation into (\ref{eq:opt_robust}), one obtains the linear approximation of the robust optimization problem:
\begin{subequations}\label{eq:opt_robust3}
	\begin{alignat}{2}
		& \underset{\pp\in \mathbb R^I}{\text{minimize}}	& \quad & J(\ppexpv) +\|\mathbf{D}\nabla_{\ppexpv}J(\ppexpv)\|_1 \qcm \\
		& \text{subject to}	& \quad & G_m(\ppexpv) +\|\mathbf{D}\nabla_{\ppexpv} G_m(\ppexpv)\|_1 \le 0 \qcm
	\end{alignat}
\end{subequations}
for $m = 1,\ldots,M$. A dual norm $||\cdot||_*$ is defined by
\begin{align}
\|\cdot\|_*: && \mathbb R^I &\rightarrow \mathbb R \nonumber\\
&&\bf{g} &\mapsto \|\bf{g}\|_* := \displaystyle\max_{\bf{g}\in\mathbb R^I,\|\dd\|\le 1} \bf{g}^\top \dd \qpt
\end{align}
In this particular case, one can use the property that the dual of $\|\mathbf{D}^{-1}\cdot\|_\infty$ is given by $\|\mathbf{D}\cdot\|_1$.

A further problem is introduced by the fact that the norms are not differentiable, which leads to a non-smooth optimization problem. A differentiable problem is obtained by introducing $M+1$ \emph{slack} variables $\slack_0,\ldots,\slack_M$ and reformulate~(\ref{eq:opt_robust3}) as
\begin{subequations}\label{eq:opt_robust_lin}
	\begin{alignat}{2}
		\label{eq:cost_smooth_robust_lin}
		& \underset{\ppexpv\in R^I, \slack_0,\ldots,\slack_M \in \mathbb R^n}{\text{minimize}}	& \quad & J(\ppexpv) +\mathbb V^\top \slack_0 \qcm \\
		\label{eq:constraints_robust_lin}
		& \text{subject to}	& \quad & G_m(\ppexpv) +\mathbb V^\top \slack_m \le 0 \qcm \\
		\label{eq:constraints_robust_lin_slack}
		& \text{}	& \quad & -\slack_0 \le \mathbf{D}\nabla_{\ppexpv} J(\ppexpv)\le \slack_0 \qcm\\
		& \text{}	& \quad & -\slack_m \le \mathbf{D}\nabla_{\ppexpv} G_m(\ppexpv)\le \slack_m \qcm 
	\end{alignat}
\end{subequations}
where $m=1,\ldots,M$ and $\mathbb V =[1,\ldots,1]^\top\in\mathbb R^I$. This optimization problem can now be efficiently solved numerically. Additionally to the quantities introduced in the previous section, now also second order sensitivities with respect to the design parameters are required. This approach can be generalized to use a quadratic approximation with respect to $\pp$ as worked out in \cite{Lass_2017aa}.

\section{Finite-Element Model}
\label{sec:model}

The behavior of the devices under consideration is determined by magnetic field phenomena and is simulated using a FE model.

\subsection{Magnetoquasistatic Formulation}
\label{subsec:formulation}

The magnetoquasistatic (MQS) subset of Maxwell's equations is considered. The design parameters $\pp$ influence the material distribution represented by the reluctivity $\nu(\pp)$ and the conductivity $\sigma(\pp)$, as well as the excitations, represented by the applied current density $\mvc{J}_{\rm src}(\pp)$ in current carrying conductors and the magnetizing field strength $\mvc{H}_{\rm m}(\pp)$ of the present permanent magnets. The MQS formulation in terms of the magnetic vector potential $\MVP(\pp)$ reads
\begin{equation}\label{eq:mqscont}
    \mvc{\nabla}\times\left(\nu(\pp)\mvc{\nabla}\times\MVP(\pp)\right)
    +\sigma(\pp)\pdt{\MVP(\pp)}
    =\mvc{J}_{\rm src}(\pp)-\mvc{\nabla}\times\mvc{H}_{\rm m}(\pp) \qcm
\end{equation}
and is complemented with adequate boundary conditions. Eq.~\ref{eq:mqscont} encompasses the case of linear, nonlinear and remanent magnetic materials expressed by
\begin{align}
	\mvc{H}(\pp) &= \nu(\pp)\mvc{B}(\pp) \qcm\\
	\mvc{H}(\pp) &= \nu(\pp,\mvc{B}(\pp))\mvc{B}(\pp) \qcm\\
	\mvc{H}(\pp) &= \mvc{H}_{\rm m}(\pp)+\nu(\pp)\mvc{B}(\pp) \qcm
\end{align}
respectively. $\mvc{H}(\pp)$ and $\mvc{B}(\pp)=\nabla\times\MVP(\pp)$ are the magnetic field strength and magnetic flux density. In the nonlinear setting, the formulation is treated by the Newton method, which is equivalent to using a linearized material relation $\mvc{H}(\pp)=\mvc{H}_{\rm m}\upk(\pp)+\overline{\overline{\nu}}\upk(\pp)\mvc{B}(\pp)$ and updating the tensorial differential permeability $\overline{\overline{\nu}}\upk(\pp)$ and the magnetizing field strength $\mvc{H}_{\rm m}\upk(\pp)$ between the successive Newton steps $k$ \cite{Koch_2008ad}.

\subsection{Finite-element discretization}
\label{subsec:discretization}

The magnetic vector potential is discretized by lowest-order N\'ed\'elec edge shape functions $\mvc{w}_j(x,y,z)$, i.e.,
\begin{equation}
	\MVP(\pp)\approx\sum_{j=1}^{N_{\rm dof}} a_j(\pp) \mvc{w}_j(x,y,z) \qcm
\end{equation}
where $a_j(\pp)$ are the degrees of freedom and $N_{\rm dof}$ is the number of degrees of freedom. In the 3D case, the shape functions are associated with the edges of a tetrahedral mesh. In the 2D cartesian case, the edge shape functions are aligned with the $z$-axis and are constructed from the nodal shape functions $N_j(x,y)$ associated with the nodes of a 2D mesh, i.e.,
\begin{equation}
	\mvc{w}_j(x,y) =\frac{N_j(x,y)}{l_z}\mvc{e}_z \qcm
\end{equation}
where $l_z$ is the length of the device in $z$-direction. In both cases, the discretization procedure leads to the system of equations
\begin{equation}
\label{eq:sys_eq}
\mathbf{K}_\nu(\pp)\mathbf{a}(\pp) +\mathbf{M}_\sigma(\pp)\fdt{\mathbf{a}(\pp)}=\mathbf{j}_{\rm src}(\pp)+\mathbf{j}_{\rm m}(\pp) \qcm
\end{equation} 
where 
\begin{align}
\label{eq:Knu} K_{\nu,i,j}(\pp) &= \int_{V_D} \nu(\pp) \mvc{\nabla}\times\mvc{w}_j\cdot\mvc{\nabla}\times\mvc{w}_i \udV \qsc\\
\label{eq:Msigma} M_{\sigma,i,j}(\pp) &= \int_{V_D} \sigma(\pp) \mvc{w}_j\cdot\mvc{w}_i \udV \qsc\\
\label{eq:jsrc} j_{\text{src},i}(\pp) &= \int_{V_D}\mvc{J}_{\rm src}(\pp)\cdot \mvc{w}_i \udV \qsc\\
\label{eq:jm} j_{\text{m},i}(\pp) &= -\int_{V_D}\mvc{H}_{\rm m}(\pp) \cdot \mvc{\nabla}\times\mvc{w}_i \udV \qcm
\end{align}
and where $V_D$ is the computational domain \cite{Monk_2003aa}. In the 2D case, $V_D=S_D \times[0,l_z]$ where $S_D$ is the cross section of the device. Eq.~\ref{eq:sys_eq} is further discretized in time by, e.g., an implicit Runge-Kutta method, linearized by the Newton-Raphson method and solved by a solution method for large sparse systems of equations \cite{Clemens_2005aa,Koch_2008ad}.

\subsection{Geometry Parametrization}
In the following, designs will be optimized with respect to geometric parameters. At first sight, the changing geometry necessitates the reconstruction of the computational mesh. This would, however, lead to unacceptably high computation times. Moreover, the unavoidable changes in mesh topology would introduce numerical noise which could mask the true sensitivity of the quantities of interest on the geometric parameters. Two different types of parametrizations are presented in the following. \emph{Affine decomposition} (see e.g. \cite{Rozza_2008aa}) is particularly appealing in the context of model order reduction and well-suited for parallelization. However, curved geometries cannot be represented exactly and additional approximation errors occur in this case. This is not the case for the second parametrization which is based on the well-established concept of \emph{design elements} \cite{Braibant_1984aa} in combination with Non-Uniform Rational B-Splines (NURBS). Here, the mapping will not be affine and more effort is needed for the update of the FE matrices and vectors. Another drawback is the difficulty in assuring the mesh quality during optimization. Yet, good results can be obtained for many shape optimization problems by one of the two methods, with moderate implementation effort. It should also be mentioned that non-parametric approaches to shape optimization \cite{Delfour_2011aa} present a viable alternative and have already been applied for electric machines \cite{Gangl_2016aa}. There, however, advanced techniques for both derivation and implementation are needed.

The geometry is decomposed in a domain $V_D^0$ that is unaffected from the geometric parameters and domains $V_D^\ell$, $\ell=1,\ldots,L$ subject to geometric changes. The FE matrices $\mathbf{K}_\nu(\pp)$ and $\mathbf{M}_\sigma(\pp)$ and vectors $\mathbf{j}_{\rm src}(\pp)$ and $\mathbf{j}_{\rm m}(\pp)$ can be partitioned accordingly, e.g.,
\begin{equation}
\label{eq:sys_eq_aff}
\mathbf{K}_\nu(\pp) =\mathbf{K}_\nu^{0}+\sum_{\ell=1}^L \mathbf{K}_\nu^{\ell}(\pp) \qcm
\end{equation}
and similarly for $\mathbf{M}_\sigma(\pp)$, $\mathbf{j}_\text{src}(\pp)$ and $\mathbf{j}_\text{m}(\pp)$. Reference geometries $\hat{V}_D^\ell, l=1,\ldots,L$ and $\hat{V}_D^0 = V_D^0$ are defined, as well as a map from $\hat{V}_D^\ell$ to $V_D^\ell$, given by $f^\ell:\hat{\mvc{r}}\rightarrow\mvc{r}=f^\ell(\hat{\mvc{r}})$, which depends on the geometric parameters $\pp$.

\subsubsection{Affine Decomposition}
\label{subsec:affine_decomposition}

In the case of affine decomposition, the domains $V_D^\ell$, $\ell=1,\ldots,L$ are triangles or tetrahedra. Hence, the maps are affine, and referred to as $f^\ell_{\rm app}$. These transformations shift the corners of the mesh, while preserving straight edges.

A key advantage of affine decomposition is that the Jacobian of the map,
\begin{equation}
	J_{\rm aff}^\ell(\pp)=\left[\begin{array}{ccc}
		\pd{x}{\hat{x}} & \pd{x}{\hat{y}} & \pd{x}{\hat{z}} \\
		\pd{y}{\hat{x}} & \pd{y}{\hat{y}} & \pd{y}{\hat{z}} \\
		\pd{z}{\hat{x}} & \pd{z}{\hat{y}} & \pd{z}{\hat{z}}
	\end{array}\right] 
\end{equation}
is constant on each subdomain $V_D^\ell$. In the integrations in (\ref{eq:Knu})-(\ref{eq:jm}), the volume integrations now have to be carried out according to $\udV=\vartheta_0^\ell(\pp)\ud{\hat{V}}$, where $\vartheta_0^\ell(\pp)=|J_{\rm app}^\ell(\pp)|$ denotes the determinant of the Jacobian. Hence,
\begin{align}
	\mathbf{M}_\sigma(\pp) &=\vartheta_0^\ell(\pp)\hat{\mathbf{M}}_\sigma \qsc\\
	\mathbf{j}_\text{src}(\pp) &=\vartheta_0^\ell(\pp)\hat{\mathbf{j}}_\text{src} \qcm
\end{align}
where $\hat{\mathbf{M}}_\sigma$ and $\hat{\mathbf{j}}_\text{src}$ are assembled for the reference geometry only once. Additionally, the affine maps affect the curl operators in (\ref{eq:Knu}) and (\ref{eq:jm}). A bit of calculation is needed to work out the transformed curl operators and the scalar products component-wise. For the 2D cartesian case, the results are
\begin{align}
\mathbf{K}_\nu^{\ell}(\pp) &= \vartheta^{\ell}_1(\pp) \hat{\mathbf{K}}^{\ell}_{\nu,xx}+\vartheta^{\ell}_2(\pp)\hat{\mathbf{K}}^{\ell}_{\nu,yy}\nonumber\\
 & \label{eq:comp_not}\quad +\vartheta^{\ell}_3(\pp) \hat{\mathbf{K}}^{\ell}_{\nu,xy}+\vartheta^{\ell}_4(\pp)\hat{\mathbf{K}}^{\ell}_{\nu,yx} \qsc\\
\mathbf{j}_\text{m}^{\ell}(\pp) &= \vartheta^\ell_5(\pp)\hat{\mathbf{j}}_{\text{m},x}^\ell+\vartheta^\ell_6(\pp)\hat{\mathbf{j}}_{\text{m},y}^\ell \qcm
\end{align}
where the matrix factors $\hat{\mathbf{K}}^{\ell}_{\nu,xx}$, $\hat{\mathbf{K}}^{\ell}_{\nu,yy}$, $\hat{\mathbf{K}}^{\ell}_{\nu,xy}$, $\hat{\mathbf{K}}^{\ell}_{\nu,yx}$, and the vector factors $\hat{\mathbf{j}}_{\text{m},x}^\ell$ and $\hat{\mathbf{j}}_{\text{m},y}$ are assembled for the reference geometry in advance. Hence, the assembly of new FE matrices and vectors can be avoided during the optimization procedure. The functions $\vartheta_q^\ell(\pp)$ are simple scalar functions in terms of the design parameters and are evaluated for each model instantiation.

\subsubsection{Design Element Approach}
\label{subsec:design_element}

NURBS are a very general way to represent geometries and widely used in CAD systems. Therefore, it seems natural to use the control points (and weights) of NURBS curves as design parameters \cite{Braibant_1984aa,Ryu_2005aa}. This approach has received considerable attention in recent years as new approaches, incorporating NURBS geometries into FE analysis, have emerged. Isogeometric analysis \cite{Hughes_2005aa} and the NURBS-enhanced FE method \cite{Sevilla_2008aa} are important examples. Here, only NURBS are used for the geometry parametrization. A triangular (tetrahedral) mesh is generated once and deformed using the well-established concept of design elements \cite{Braibant_1984aa,Imam_1982aa}. 

In the following, for simplicity, the two-dimensional case is considered solely. A generic NURBS curve of degree $p$ is given as 
\begin{equation}
\mvc{C}(\hat{x}) = \sum_{i} R_i^p(\hat{x}) \mvc{P}_{i} \qcm
\end{equation}
where $\mvc{P}_{i}$ refers to a control point and the rational spline $R_i^p$ is defined in terms of B-splines $N_i^p$ and weights $w_i$ as 
\begin{equation}
R_i^p(\hat{x}) = \frac{N_i^p(\hat{x}) w_i}{\sum_j N_j^p(\hat{x}) w_j} \qpt
\end{equation}
In total, $L$ design elements are considered, each of which is represented by two NURBS curves $\mvc{C}_1^\ell$ and $\mvc{C}_2^\ell$. More precisely, a design element is defined by a map $f_{\rm de}^\ell:\hat{V}_D^\ell = [0,1]^2\rightarrow V_D^\ell$ given as
\begin{equation}
f_{\rm de}^\ell(\hat{x},\hat{y}) = \mvc{C}_1^\ell(\hat{x}) \hat{y} + \mvc{C}_2^\ell(\hat{x}) (1 - \hat{y}) \qpt
\end{equation}
Hence, design elements are given as Cartesian products of NURBS curves, whereas the affine decomposition may result in unstructured representations. For each node $(x_i,y_i)$ in $V_D^\ell$, its position in the reference domain $[0,1]^2$ is computed in advance by solving 
\begin{equation}
(\hat{x}_i,\hat{y}_i) \in {\rm argmin}_{(\hat{x},\hat{y})} |f_{\rm de}^\ell(\hat{x},\hat{y}) - (x_i,y_i)| \qcm
\end{equation}
e.g., with the Newton-Raphson method. Then, the mesh can be easily deformed by applying the parameter-dependent map $f_{\rm de}^\ell$ to all nodes $(\hat{x}_i,\hat{y}_i)$.  

The transformation of the FE matrices and vectors is more involved compared to the affine decomposition described in Section~\ref{subsec:affine_decomposition}. Each entry of the mass matrix is transformed as 
\begin{equation}
\mathbf{M}_{\sigma,i,j}(\pp) = \int_{\hat{V}_D} \hat{\sigma} \hat{\mvc{w}}_j\cdot \hat{\mvc{w}}_i |J_{\rm de}^\ell(\pp)| \ud{\hat{V}},
\end{equation}
where it is important to emphasize that $|J_{\rm de}^\ell(\pp)|$ is not constant on each design element. A similar expression is obtained for $\mathbf{j}_\text{src}(\pp)$, whereas the conforming transformation of the curl operator yields
\begin{align} \label{eq:Knutransformed}
\mathbf{K}_{\nu,i,j}(\pp) &= \int_{\hat{V}_D} \frac{\hat{\nu}}{|J_{\rm de}^{\ell}(\pp)|} J_{\rm de}^{\ell}(\pp) \mvc{\nabla} \times \hat{\mvc{w}_j}\cdot J_{\rm de}^\ell(\pp) \mvc{\nabla} \times \hat{\mvc{w}_i}  \ud{\hat{V}}, \\
\mathbf{j}_{\text{m},i}^{\ell}(\pp) &= \int_{\hat{V}_D}\hat{\mvc{H}}_{\rm m}(\pp) \cdot J_{\rm de}^{\ell}(\pp) \mvc{\nabla}\times \hat{\mvc{w}_i} \ud{\hat{V}}. \label{eq:jtransformed}
\end{align}
In \eqref{eq:Knutransformed} and \eqref{eq:jtransformed}, the dependence of the integration domain on the geometry changes was eliminated. Because the Jacobian $J_{\rm de}^{\ell}(\pp)$ can be expressed as a function of the geometry parameters $p_i$, the analytical derivative of the system matrix and of the right hand side with respect to the geometry parameters can be determined.

\subsection{Sensitivities}

After differentiating the FE system, a new linear system for the derivatives of the degrees of freedom with respect to the geometry parameters is obtained:
\begin{equation}
	\label{eq:sens}
	\mathbf{K}_\nu\mathbf{s}_i = \pd{}{p_i}\left(\mathbf{j}_\text{src}+\mathbf{j}_\text{m}\right)
	-\pd{\mathbf{K}_\nu}{p_i}\mathbf{a} \qcm\quad\text{for }i = 1,\ldots,I \qcm
\end{equation}
where $\mathbf{s}_i(\pp)=\pd{{\bf a}({\pp})}{p_i}$ are the sensitivities of the FE solution. To calculate $\mathbf{s}_i$, $I$ equations of the form \eqref{eq:sens} have to be solved. In the case of affine decomposition, derivatives of $\mathbf{K}_\nu$ are easily calculated from \eqref{eq:sys_eq_aff} and \eqref{eq:comp_not} using expressions for $\pd{\vartheta^\ell(\pp)}{p_i}$ which are known analytically as derivatives of the functions 
$\vartheta^{\ell}(\pp)$. The expressions become more involved when NURBS are involved, yet closed form formulas also exist in this case.

The optimization algorithm requires the derivatives of the objective function with respect to each of the design parameters. Often, the objective function does not explicitly depend on the design parameters, i.e., $J(\pp)=\tilde{J}(\bf a(\pp))$. In this case, the derivatives are given as
\begin{equation}
	\label{eq:sensder}
	\pd{J(\pp)}{p_i} = \nabla_{\bf a}\tilde{J}({\bf a(\pp)}) \cdot {\bf s}_i(\pp) \text{ , for }i = 1,\ldots,I.
\end{equation}
For a large number of parameters, an adjoint method should be used instead \cite{Troltzsch_2009aa}.

\section{Example~1: Die Press Mold}
\label{sec:example1}

As a first example, a die press mold for radially magnetizing a segment of sintered magnetic powder (SMP) is considered \cite{Takahashi_1996aa}. This problem has been proposed as Testing-Electromagnetic-Analysis-Methods (TEAM) benchmark problem~25 \cite{TEAM25_1997aa} and has been used in numerous papers for comparing optimization algorithms. The vast majority of these publications apply and compare \emph{stochastic} optimization methods \cite{Lei_2015aa,Sonoda_2007aa}, possibly combined with surrogate models \cite{Canova_2003aa}, uncertainty quantification \cite{Nishida_2008aa}, multi-objective optimization or a combination of them \cite{Lebensztajn_2004aa}. Only a few papers (e.g. \cite{Alotto_2001aa} and \cite{Berkani_2013aa}) choose deterministic methods, again possibly combined with surrogate models \cite{Hemker_2007aa}, uncertainty quantification \cite{Takahashi_2001aa} or multi-objective optimization. This paper addresses one of the main drawbacks of deterministic methods, i.e., the consideration of geometric parameters. For this example, the design element approach is used. 

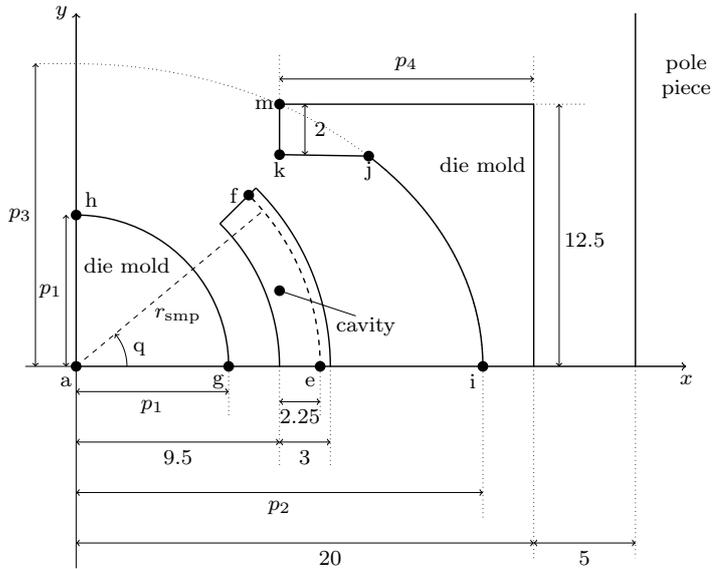
\begin{figure}[tb]
	\centering\small
	\resizebox{.8\columnwidth}{!}{\begin{tikzpicture}[scale=1, every node/.style={scale=1.5}]
\XYcoord{12}{7};
\draw (0,0) -- (-1,0);
\draw (0,0) -- (0, -4);

\node [fill,circle, scale=0.5] at (0,0) (aa) {};
\node at (-0.2, -0.3) (a) {a};
\draw [thick,domain=0:90] plot ({0+3*cos(\x)}, {0+3*sin(\x)});
\node [fill,circle, scale=0.5] at ({0+3*cos(0)}, {0+3*sin(0)}) (gg) {};
\node at ({-0.2+3*cos(0)}, {-0.3+3*sin(0)}) (g) {g};
\node [fill,circle, scale=0.5] at ({0+3*cos(90)}, {0+3*sin(90)}) (hh) {};
\node at ({+0.3+3*cos(90)}, {+0.3+3*sin(90)}) (h) {h};

\draw [thick,domain=0:45] plot ({4*cos(\x)}, {4*sin(\x)});
\draw [thick,domain=0:45] plot ({5*cos(\x)}, {5*sin(\x)});
\draw [thick] ({4*cos(45)}, {4*sin(45)}) -- ({5*cos(45)}, {5*sin(45)});
\draw [thick, dashed, domain=0:45] plot ({4.8*cos(\x)}, {4.8*sin(\x)});
\node [fill,circle, scale=0.5] at ({4.8*cos(0)}, {4.8*sin(0)}) (ee) {};
\node at ({-0.2+4.8*cos(0)}, {-0.3+4.8*sin(0)}) (e) {e};
\node [fill,circle, scale=0.5] at ({4.8*cos(45)}, {4.8*sin(45)}) (ff) {};
\node at ({-0.3+4.8*cos(45)}, {4.8*sin(45)}) (f) {f};

\draw [dotted, domain=0:90] plot ({8*cos(\x)}, {6*sin(\x)});
\node [fill,circle, scale=0.5] at ({8*cos(0)}, {6*sin(0)}) (ii) {};
\node at ({-0.2+8*cos(0)}, {-0.3+6*sin(0)}) (i) {i}; 
\node [fill,circle, scale=0.5] at ({8*cos(60)}, {6*sin(60)}) (mm) {};
\node at ({-0.3+8*cos(60)}, {6*sin(60)}) (m) {m};
\node [fill,circle, scale=0.5] at ({8*cos(60)}, {-1+6*sin(60)}) (kk) {};
\node at ({8*cos(60)}, {-1.3+6*sin(60)}) (k) {k};
\draw [thick] ({8*cos(60)}, {6*sin(60)}) -- ({8*cos(60)}, {-1+6*sin(60)});
\draw [thick] ({8*cos(60)}, {-1+6*sin(60)}) -- ({8*cos(44)}, {6*sin(44)});
\node [fill,circle, scale=0.5] at ({8*cos(44)}, {6*sin(44)}) (jj) {};
\node at ({8*cos(44)}, {-0.3+6*sin(44)}) (j) {j};
\draw [thick, domain=0:44] plot ({8*cos(\x)}, {6*sin(\x)});
\draw [thick] ({1+8*cos(0)}, {6*sin(60)}) -- ({8*cos(60)}, {6*sin(60)});
\draw [thick] ({1+8*cos(0)}, {6*sin(60)}) -- ({1+8*cos(0)}, {6*sin(0)});

\draw [thick] (11,0) -- (11,7);

\draw [dashed] (0,0) -- ({4.8*cos(40)}, {4.8*sin(40)});
\node at (2, 1) (Rsmp) {$r_{\text{smp}}$};
\draw [->, domain=0:40] plot ({cos(\x)}, {sin(\x)});
\node at ({0.3+cos(20)}, {sin(20)}) (q) {q};

\node at (12, 6) (pp) {pole};
\node at (12, 5.5) (pp) {piece};
\node at (1, 2) (dm1) {die mold};
\node at (8, 4) (dm2) {die mold};
\node[fill,circle, scale=0.5] at (4, 1.5) (cv) {};
\draw (4, 1.5) -- (5.5, 1);
\node at (5.7, 0.8) (cv2) {cavity};

\draw [<->] ($(aa) - (0.2, 0)$) -- ($(hh) - (0.2, 0)$);
\node at (-0.5, 1.5) (R1y) {$p_1$};
\draw [<->] ($(aa) - (0.8, 0)$) -- ({-0.8+8*cos(90)}, {6*sin(90)});
\node at (-1.1, 3) (L3) {$p_3$};
\draw [dotted] ({-0.8+8*cos(90)}, {6*sin(90)}) -- ({8*cos(90)}, {6*sin(90)});
\draw [<->] ($(aa) - (0, 0.5)$) -- ($(gg) - (0, 0.5)$);
\node at (1.5, -0.8) (R1x) {$p_1$};
\draw [dotted] (gg) -- ($(gg) - (0, 1)$);
\draw [<->] ($(aa) - (0, 1.5)$) -- ($({4*cos(0)}, {4*sin(0)}) - (0, 1.5)$);
\node at (2, -1.8) (9.5) {$9.5$};
\draw [dotted] ({4*cos(0)}, {4*sin(0)}) -- ($({4*cos(0)}, {4*sin(0)}) - (0, 2)$);
\draw [<->] ($({4*cos(0)}, {4*sin(0)}) - (0,1.5)$) -- ($({5*cos(0)}, {5*sin(0)}) - (0,1.5)$);
\draw [dotted] ({5*cos(0)}, {5*sin(0)}) -- ($({5*cos(0)}, {5*sin(0)}) - (0, 2)$);
\node at (4.5, -1.8) (3) {$3$};
\draw [<->] ($({4*cos(0)}, {4*sin(0)}) - (0, 0.7)$) -- ($(ee) - (0,0.7)$);
\draw [dotted] (ee) -- ($(ee) - (0, 1)$);
\node at (4.4, -1) (2.25) {$2.25$};
\draw [<->] ($(aa) - (0, 2.5)$) -- ($(ii) - (0, 2.5)$);
\node at (4, -2.8) (L2) {$p_2$};
\draw [dotted] (ii) -- ($(ii) - (0, 3)$); 
\draw [<->] ($(aa) - (0, 3.5)$) -- ($({1+8*cos(0)}, {6*sin(0)}) - (0, 3.5)$);
\draw [dotted] ($({1+8*cos(0)}, {6*sin(0)}) - (0, 3.8)$) -- ({1+8*cos(0)}, {6*sin(0)});
\node at (5, -3.8) (20) {$20$};
\draw [dotted] (11,0) -- (11,-3.8);
\draw [<->]  ($({1+8*cos(0)}, {6*sin(0)})-(0,3.5)$) -- (11,-3.5);
\node at (10, -3.8) (5) {$5$};
\draw [<->] ($({1+8*cos(0)}, {6*sin(0)}) + (0.5, 0)$) -- ($({1+8*cos(0)}, {6*sin(60)}) + (0.5, 0)$);
\draw [dotted] ($({1+8*cos(0)}, {6*sin(60)}) + (1, 0)$) -- ({1+8*cos(0)}, {6*sin(60)});
\node at (10, 2.5) (12.5) {$12.5$};
\draw [<->] ($(mm) + (0, 0.5)$) -- ($({1+8*cos(0)}, {6*sin(60)}) + (0, 0.5)$);
\draw [dotted] (mm) -- ($(mm) + (0, 1)$); 
\draw [dotted] ({1+8*cos(0)}, {6*sin(60)}) -- ($({1+8*cos(0)}, {6*sin(60)}) + (0, 1)$);
\node at (6.5, 6) (L4) {$p_4$};
\draw [<->] ($(mm) + (0.5, 0)$) -- ($(kk) + (0.5, 0)$);
\node at ($(kk) + (0.8, 0.5)$) (2) {$2$};
\end{tikzpicture}}
	\caption{TEAM Problem 25: Cross section of the inner part of the die press showing the SMP ring, the inner yoke and the outer yoke (all measures in mm). A horizontal magnetic flux is exerted on the configuration by an outer magnetic circuit (not shown).}
	\label{fig:TEAM25}
\end{figure}

The SMP segment is arranged between a cylindrical inner pole and a more generally shaped outer pole (Fig.~\ref{fig:TEAM25}). The original TEAM-25 problem considers an outer pole with an elliptical inner surface. Here, the inner surface is described by a spline. This is motivated by the fact that splines are currently the basic building block for mechanical processing. The considered design parameters are then chosen to be
\begin{align}
	p_1 &: \textrm{radius of the inner yoke} \nonumber \qsc\\
	p_2,p_3 &: \textrm{semiaxis of ellipse between points i and j} \nonumber \qsc\\
	p_4 &: \textrm{$x$-coordinate of points m and k} \nonumber \qpt
\end{align}
Both the circle and the ellipse are exactly represented by NURBS curves. The relation between the geometric parameters and the NURBS control points is given in the appendix.

The optimization aims at a homogeneous, radially oriented magnetic flux density of $B_{\rm goal}=0.35$~T inside the SMP segment. The objective function $J(\pp)$ is defined as the mean squared error between the simulated magnetic field and the goal at $9$ sample points equidistantly distributed along the arc with radius $r_\text{smp}$ between points $e=(r_\text{smp},0)$ and $f=(r_\text{smp}\cos\varphi_f,r_\text{smp}\sin\varphi_f)$, i.e.,
\begin{equation} \label{eq:J}
	J(\pp) = \sum_{k=1}^9 \| \mvc{B}(r_{\rm smp}\cos\varphi_k,r_{\rm smp}\sin\varphi_k;\pp)-B_{\rm goal}\mvc{e}_k\|_2^2 \qcm
\end{equation}
where $\varphi_k=\varphi_f \frac{k-1}{8}$ and $\mvc{e}_k=(\cos\varphi_k,\sin\varphi_k)$.
The optimization problem yields:
\begin{subequations}
	\begin{alignat}{2}
		& \underset{\pp}{\text{minimize}}	& \quad &  J(\pp) \qcm \\
		& \text{subject to}	& \quad & \pp \in \mathcal{F}\qcm
	\end{alignat}
\end{subequations}
where the admissible set is defined as:
\begin{equation*}
\mathcal{F}=[5.1,9] \times [16,18] \times [14.5,16] \times [9.5,13] \text{ mm}.
\end{equation*} 
For the gradient, the derivatives of $J$ with respect to the geometry parameters $p_i$ are needed. Before applying the chain rule on \eqref{eq:J}, the derivatives of the degrees of freedom with respect to the geometry parameters $\partial_{p_i} \mathbf{a}$ are calculated described in Section~\ref{subsec:design_element}. 

The performance of a standard algorithm for Particle Swarm Optimization (PSO), of the Sequential Quadratic Programming (SQP) method implemented in MATLAB\textsuperscript{\textregistered}'s {\tt fmincon} function \cite{Powell_1978aa} and of an own implementation of SQP are compared in Table~\ref{tb:results}. Both SQP implementations use the analytical gradients, the BFGS formula for updating the Hessian and a sufficient decrease condition in a merit function. For the PSO, a set of $40$ particles is considered and the implementation is multi-threaded, while the gradient-based methods are single-thread implementations. The termination criterion for the PSO algorithm is the number of stall iterations, which was set to~5. The PSO actually finds the optimum after 2~iterations. This is because the optimum is at a vertex of the box-shaped domain and all the particles leaving the admissible region are projected onto the boundary. All three methods converge to the same optimum. The deterministic algorithms are by substantially faster than PSO, even though PSO exploits parallelization. On the same machine, an evaluation of the objective function $J(\pp)$ is performed in 1.65~s, an analytical evaluation of the gradient $\nabla J(\pp)$ in 4.69~s and a numerical evaluation of the gradient $\nabla_{\text{num}} J(\pp)$ using a forward difference quotient in 7.48~s. All tests were done on a 64~GB RAM Intel\textsuperscript{\textregistered} Xeon\textsuperscript{\textregistered} E5-2630 v4 machine.

\begin{table*}[tb]
	\centering\small
	\caption{Results from the optimization of the die press mold with particle swarm optimization (PSO), trust region (TR) (with MATLAB\textsuperscript{\textregistered}'s {\tt fmincon}) and an own implementation of sequential quadratic programming (SQP) combined with the design element approach.}
	\label{tb:results}
	\begin{tabular}{ccccccc}
			\hline
			method & minimizer $\hat{\pp}_{\rm min}$ & minimum & iteration &  \multicolumn{2}{c}{function calls} & time\\
			& (in mm) & (in $\text{T}^2$) & count & $\mathrm{f}()$ & $\nabla \mathrm{f}()$ & (in s)\\
			\hline
			PSO           & $\begin{pmatrix} 5.1000\\16.0000 \\16.0000\\9.5000\end{pmatrix}$ &  1.413498 & 7 & 280 & N/A & 56.63\\			           
			\hline	 		
				\begin{tabular}{@{}c@{}}SQP \\  (\texttt{fmincon})\end{tabular} & $\begin{pmatrix} 5.1000\\16.0000 \\16.0000\\9.5000\end{pmatrix}$ & 1.413498 & 4 & 7 & 7 & 31.61\\
			\hline
			
			\begin{tabular}{@{}c@{}}SQP \\ (own\\implementation)\end{tabular}      & $\begin{pmatrix} 5.1000\\16.0000 \\16.0000\\9.4999\end{pmatrix}$ & 1.413498 & 2 & 3 & 2 & 12.84\\
			
			\hline
		\end{tabular}
	\end{table*} 

\section{Example~2: Permanent-Magnet Synchronous Machine (PMSM)}
\label{sec:example2}

\subsection{Design parameters}

\begin{figure}[tb]
	\centering
	\def\svgwidth{0.99\columnwidth}
	\input{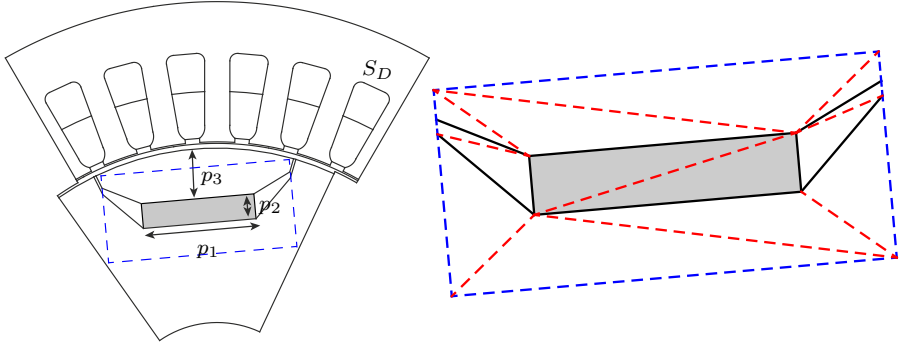}
	\caption{Cross-section of one pole of the machine with the magnet depicted in gray and the region of the affine decomposition indicated by the dashed box. On the right hand side, the triangulation into $L $ subdomains is shown by the dashed-dotted lines. Figure adapted from~\cite{Bontinck_2017ab}.}
	\label{fig:aff_dec}
\end{figure}

The second example is a $3$-phase $6$-pole permanent-magnet (PM) synchronous machine (PMSM) borrowed from \cite{Pahner_1998aa} (Fig.~\ref{fig:aff_dec}) and already studied as an optimization example in~\cite{Bontinck_2017ab}. The stator features two slots per pole and per phase with a conventional distributed double-layer winding. The rotor contains a buried rare-earth magnet. The yoke parts are laminated. The design parameters are
\begin{align}
	p_1 &: \textrm{width of the PM} \qsc \nonumber\\
	p_2 &: \textrm{thickness of the PM} \qsc \nonumber\\
	p_3 &: \textrm{distance from the PM to the rotor surface} \qpt \nonumber
\end{align}

\subsection{Objective function}

The optimization goal is to minimize the size $S_{\rm pm}=p_1p_2$ of PM material while preserving a prescribed electromotive force $E_0$. The electromotive force (EMF) $E_0(\pp)$ is post-processed from a magnetostatic solution of a 2D FE model of the PMSM using the \emph{loading} method proposed in \cite{Rahman_1991aa}. For that purpose, the FE solution of the $z$-component of the magnetic vector potential is sampled at a circle (or in the case of a partial machine model, an arc) in the PMSM's air gap, yielding $A_z(r_{\rm ag},\varphi)\approx\hat{A}_{z,{\rm eff}}\sqrt{2}\sin(N_p\varphi-\varphi_{\rm d})$, where $N_{p}=3$ is the pole-pair number, $\hat{A}_{z,{\rm eff}}$ is the rms magnitude of fundamental harmonic component and $\varphi_{\rm d}$ is the angle of the PMSM's \emph{direct} axis. The EMF is then found from
\begin{equation}
	E_0 = 2\hat{A}_{z,{\rm eff}}\omega_{\rm syn}N_{\rm w}k_{{\rm w},1} \qcm
\end{equation}
where $\omega_{\rm syn}$ is the synchronous speed and $N_{\rm w}$ is the number of windings per phase. The winding factor is
\begin{equation}
	k_{{\rm w},\nu} =\frac{\sin\left(q\nu\frac{\alpha_{\rm el}}{2}\right)}{q\sin\left(\nu\frac{\alpha_{\rm el}}{2}\right)}
	\cdot\sin\left(\nu\frac{\pi}{2}\frac{\tau_{\rm c}}{\tau_{\rm p}}\right)
	\cdot\frac{\sin\left(\nu\frac{\varepsilon}{2}\right)}{\nu\frac{\varepsilon}{2}} \qcm
\end{equation}
where $q$ is the number of coil sides per phase belt, $\alpha_{\rm el}$ is the electric angle between two slots, $\tau_{\rm c}$ is the coil pitch, $\tau_{\rm p}$ is the pole pitch and $\varepsilon$ is the electric skew angle \cite{Mueller_2008aa,Rahman_1991aa}.

\subsection{Optimization problem}

The optimization problem reads
\begin{subequations}\label{eq:optprbpmsm}
\begin{alignat}{2}
	& \underset{\pp\in \mathbb R^3}{\text{minimize}}	& \quad &  J(\pp)=p_1 p_2 l_z \qcm \\
		& \text{subject to}	& \quad & G(\pp) =\left[\begin{array}{c}
		   \pplow_1 - p_1\\
		   \pplow_2  - p_2\\
		   \pplow_3  - p_3\\
		   p_3 - \ppup_3\\
		   p_2 + p_3 - 15\;\textrm{mm}\\
		   3p_1 - 2p_3 - 50\;\textrm{mm}\\  
		   E_\mathrm{d} - E_0(\ppexpv,{\bf a}(\ppexpv))
		   \end{array}\right] \leq 0\qpt
\end{alignat}
\end{subequations}
The first four constraints are related to the lower ($\pplow$) and upper ($\ppup$) bounds of $\pp$: $(\pplow_1,\pplow_2,\pplow_3) = (1,1,5)$~mm and $(\ppup_1, \ppup_2, \ppup_3) =(\infty, \infty, 14)$~mm. To ensure the validity of the affine decomposition, i.e., intersections are not allowed, the fifth constraint is added. The sixth constraint is a design constraint enforcing that each PM has to keep a sufficient distance to the rotor surface, especially for wide PMs. The last constraint expresses the requirement to fulfill the prescribed EMF. Since the EMF is post-processed from the FE solution, the optimization problem actually has a PDE constraint.

\subsection{Results}

The results for 5~different optimization methods are collected in Table~\ref{tab:res_OP}.
\begin{enumerate}
	\item The first optimization run is carried out with the genetic algorithm implemented in MATLAB\textsuperscript{\textregistered}.
	\item The second optimization run is carried out with MATLAB\textsuperscript{\textregistered}'s PSO implementation. To circumvent the restriction to box-shaped parameter domains, the admissible set is enforced by a penalty turn. The new objective function reads
\begin{align}
	\nonumber J_\text{pen}(\pp) &=J(\pp) +2J(\pp) \big(
		f\left(\max(p_2+p_3-15,0)\right) \\ 
	\nonumber	&\quad +f\left(\max(3p_1-2p_3-50,0)\right) \\
		&\quad +f\left(\max(g(x),0)\right) \big) \qcm
\end{align}
where $f(t)=e^{(4t^{0.1})}-1$ was chosen heuristically such that $J_\text{pen}$ grows exponentially if one of the constraints is violated. The function $J_\text{pen}$ was called 4740 times, but was organized as to only evaluate the nonlinear constraint if all other constraints were satisfied. The number of particles was set to 30, the maximum number of stall iterations to $N_\text{stall}=15$ and the function change tolerance to $10^{-6}$. The PSO characteristic constants are chosen to be $\omega_0=0.5$ and $\omega_1=\omega_2=1.49$. The algorithm took 157 iterations before termination.
\item The third optimization is carried out with an own PSO implementation, for the original objective function $J(\pp)$ and applying the nonlinear constraints directly. Here, it is assumed that the admissible set is convex such that points inside the convex hull formed by all previous points do not need to be checked. 50~particles were used. Termination was enforced after maximally $N_\text{it,max}=100$ steps or when $N_\text{stall,max}=15$ stall iterations were observed.
\item The fourth run was done with a deterministic method, relying upon FE simulations equipped with an affine decomposition of the geometry as described in subsection~\ref{subsec:affine_decomposition}.
\item The fifth run was done with a deterministic method for robust optimization, again with affine decomposition of the geometry.
\end{enumerate}
The three stochastic algorithms were run on a $\unit{64}{GB}$ RAM Intel\textsuperscript{\textregistered} Xeon\textsuperscript{\textregistered} E5-2630 v4 machine. Both deterministic algorithms were run on a $\unit{16}{GB}$ RAM Intel\textsuperscript{\textregistered} Core\textsuperscript{\texttrademark} with i7-5820K processors ($\unit{3.30}{GHz}$).

The results of all optimization procedures are compared with the values of the initial design. All routines achieve a substantial decrease of the PM size from $\unit{133}{mm^2}$ up to about $\unit{63}{mm^2}$. The price for robustness is a slightly larger size of about $\unit{77}{mm^2}$. The deterministic methods outperform the stochastic ones by two orders of magnitude. This impressively illustrates the major message of this paper stating that deterministic optimization methods accompanied by FE analysis providing gradients with respect to geometric parameters should be favored over stochastic methods, at least for the here considered class of problems.

\begin{figure*}[t]
  \begin{subfigure}{0.30\textwidth}
    \centering
    \includegraphics[width=0.90\textwidth]{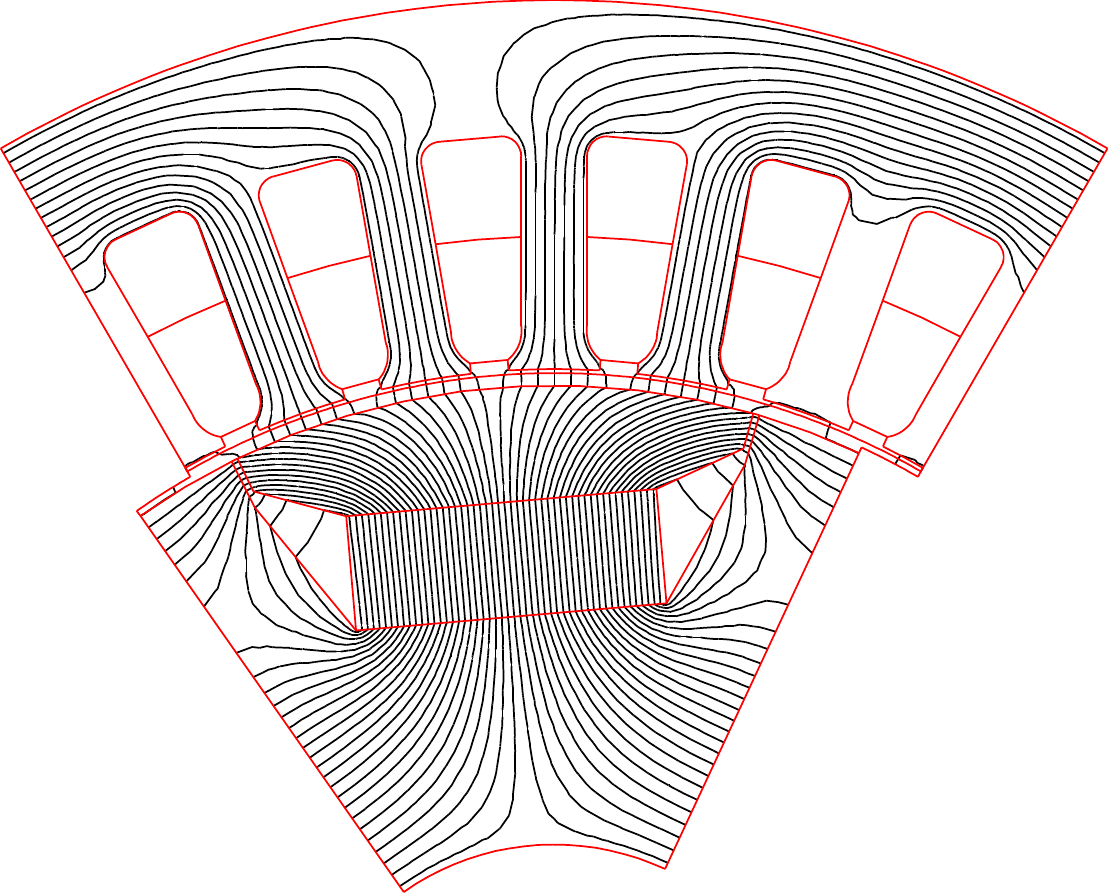}
    \caption{initial geometry}
    \label{fig:initial}
  \end{subfigure}
  \hspace{0.01\textwidth}
  \begin{subfigure}{0.30\textwidth}
    \centering
    \includegraphics[width=0.90\textwidth]{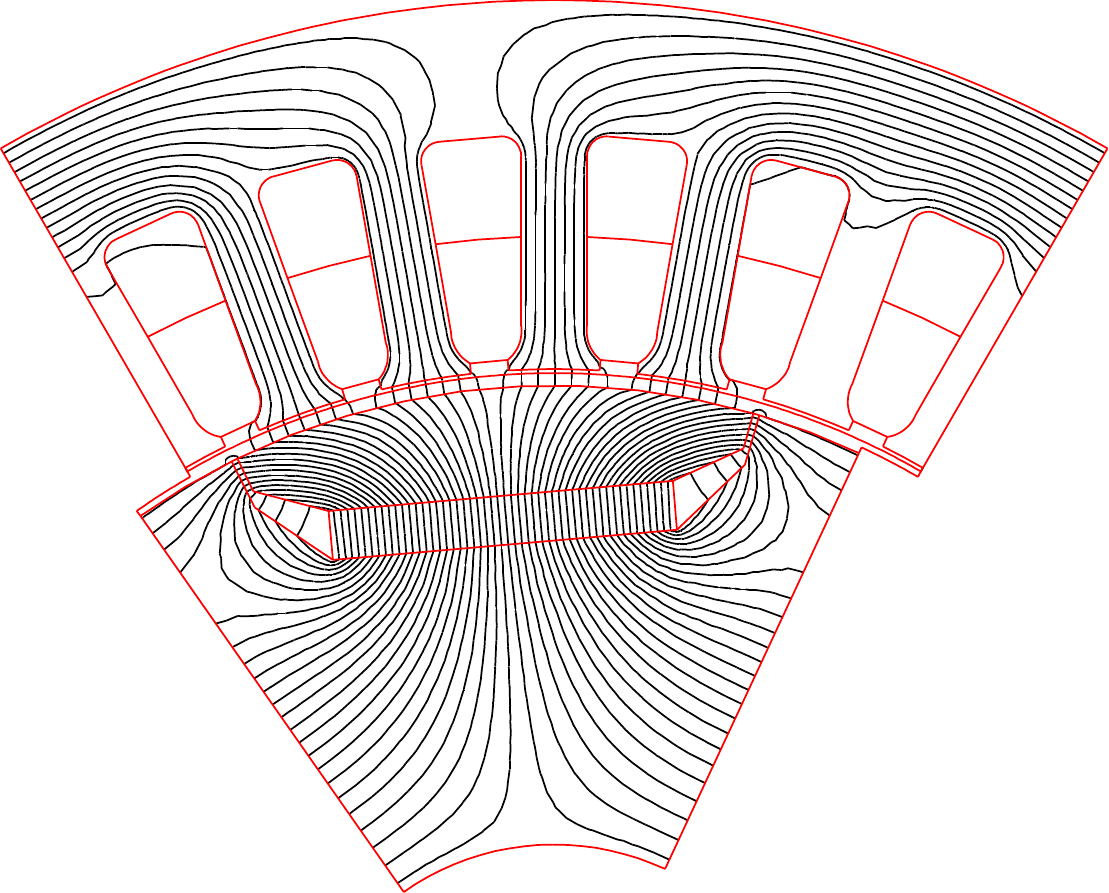}
    \caption{nominal optimum}
    \label{fig:nomopt}
  \end{subfigure} 
  \begin{subfigure}{0.30\textwidth} 
	\centering
  	\includegraphics[width=0.90\textwidth]{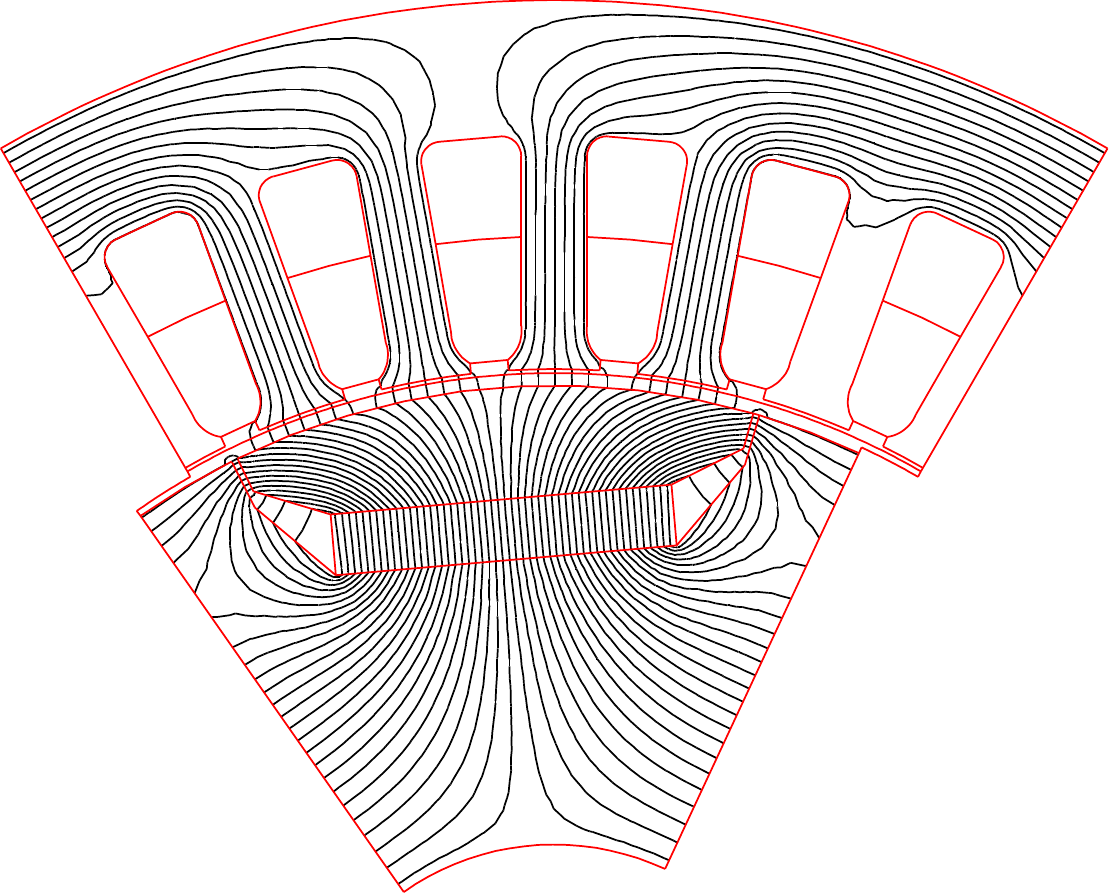}
    \caption{robust optimum}
    \label{fig:robopt}
  \end{subfigure}
  \caption{Initial and optimized geometries together with the magnetic flux distribution at no-load. Figures adapted from~\cite{Bontinck_2017ab}.}
  \label{fig:results}
\end{figure*}

\begin{table*}[tb]
\caption{Numerical results obtained for a $\dd=\unit{0.2}{mm}$~\cite{Bontinck_2017ab}.}
\label{tab:res_OP}
\begin{center}
\begin{tabular}{p{2.5cm}ccccccc}
\hline                           & $p_1$ & $p_2$ & $p_3$ & $S_{\rm pm}$      & $E_0$  & FE slv         & time   \\
                                 & (mm)  & (mm)  & (mm)  & $(\unit{}{mm^2})$ & (V)    & calls          & (s)    \\
\hline initial design            & 19.00 & 7.00  & 7.00  & 133               & 30.370 & -              & -      \\
\hline genetic algorithm         & 21.04 & 2.98  & 6.56  & 62.80             & 30.370 & $\approx 6760$ & 520.5  \\
\hline PSO with penalty term     & 20.60 & 3.09  & 5.91  & 63.71             & 30.370 & $\approx 3470$ & 267.16 \\
\hline PSO, own implementation   & 21.08 & 2.98  & 6.63  & 62.80             & 30.370 & 1765           & 217.52 \\
\hline SQP, nominal optimization & 21.07 & 2.98  & 6.61  & 62.80             & 30.370 & 34             & 2.0    \\
\hline SQP, robust optimization  & 20.88 & 3.73  & 6.82  & 77.87             & 31.086 & 48             & 5.9    \\
\hline
\end{tabular}
\end{center}
\end{table*}
 
\section{Conclusion}
\label{sec:conclusions}

Affine decomposition and design element approaches are capable of parametrizing the geometry of finite-element models such that accurate derivatives with respect to geometric parameters become available. This alleviates one of the major drawbacks of gradient-type deterministic optimization methods. For the example of a die mold press, standard sequential programming combined with the design element approach outperforms particle swarm optimization by more than a factor ten. The second example illustrates the applicability of gradient-type robust optimization combined with an affine decomposition of the geometry for a permanent-magnet synchronous machine. Supported by the substantial improvement in computational efficiency, this paper stands up for a revival of deterministic methods for numerical optimization in electrotechnical design procedures.

\section{Appendix}
\label{sec:appendix}
The dependence between the geometry parameters and the NURBS representation of the die press model is as follows. For the ellipse arc, the control points and weights are
\begin{equation*}
\mathbf{P_0}=\begin{pmatrix}
L_2\\0
\end{pmatrix} \qcm\quad
\mathbf{P_2}=\begin{pmatrix}
L_2 \cos{\alpha}\\ L_3 \sin{\alpha} 
\end{pmatrix} \qcm
\end{equation*}
\begin{equation*}
\mathbf{P_1}=\mathbf{P_2}+\begin{pmatrix}
-\lambda L_2 \sin{\alpha}\\ \lambda L_3 \cos{\alpha}
\end{pmatrix} \qcm
\end{equation*}
\begin{equation*}
w_0=w_2=1 \qcm\quad \hspace{1cm} w_1=\cos{\frac{\alpha}{2} } \qcm
\end{equation*}
where
\begin{equation*}
\alpha=\text{asin}\left(\frac{10.5mm}{p_3}\right) \qcm\quad
\lambda=\frac{-p_2+p_2\cos\alpha}{p_2\sin\alpha} \qpt
\end{equation*}
The corresponding knots are $\mathcal{K}=\{0,0,0,1,1,1\}$ and the degree of the basis functions is $p=2$.

For the circular arc, the control points are
\begin{equation*}
\mathbf{P_0}=\begin{pmatrix}
p_1\\0
\end{pmatrix} \qcm\quad
\mathbf{P_1}=\begin{pmatrix}
p_1\\p_1
\end{pmatrix} \qcm\quad
\mathbf{P_2}=\begin{pmatrix}
0\\p_1
\end{pmatrix} \qcm
\end{equation*}
with the constant weights
\begin{equation*}
w_0=w_2=1 \qcm\quad w_1=\sqrt{2}/2 \qpt
\end{equation*}
The degree of the basis functions is $p=2$. The corresponding knots are $\mathcal{K}=\{0,0,0,1,1,1\}$. The deformation of the mesh inside one design element region $V_{D}^{\ell}$ with $N_V$ vertices is given by
\begin{equation}
 (x_i,y_i)=\mvc{C}_1^\ell(\hat{x}_i;\pp) \hat{y}_i + \mvc{C}_2^\ell(\hat{x}_i;\pp) (1 - \hat{y}_i) \qcm
\end{equation}
where $(x_i,y_i)$ are the coordinates of the vertices of the deformed mesh and $(\hat{x}_i,\hat{y}_i)$ are the coordinates in the reference domain $[0,1]^2$.

\section*{Acknowledgment}
This work is supported by the German BMBF in the context of the SIMUROM project (grant nr. 05M2013) and the PASIROM project (grant nr. 05M2018), by the 'Excellence Initiative' of the German Federal and State Governments and by the Centre and Graduate School Computational Engineering at TU Darmstadt.

\end{document}